\documentclass[11pt,twoside]{amsart}

\usepackage{amsmath}
\usepackage{amsthm}
\usepackage{amsfonts, amssymb}
\usepackage{mathrsfs}
\usepackage[all,line]{xy}
\usepackage{url}
\usepackage{graphics}
\usepackage{epstopdf}

\usepackage{latexsym}
\usepackage[dvips]{graphicx}

\theoremstyle{plain}
\newtheorem{lemma}{Lemma}[section]
\newtheorem{proposition}[lemma]{Proposition}
\newtheorem{thm}[lemma]{Theorem}
\newtheorem{corollary}[lemma]{Corollary}
\theoremstyle{remark}

\newtheorem{remark}[lemma]{Remark}

\theoremstyle{definition}
\newtheorem{definition}[lemma]{Definition}
\newtheorem{ej}[lemma]{Example}

\def\u#1{U_{#1}}
\def\f#1{F_{#1}}
\def\ft#1{\hat{F}_{#1}}
\def\ut#1{\hat{U}_{#1}}
\def\ct#1{\hat{C}_{#1}}
\def\set#1{\{#1\}}
\def\k{\mathcal{K}}
\def\x{\mathcal{X}}
\def\h{\mathcal{H}}
\def\Z{\mathbb{Z}}
\def\R{\mathbb{R}}
\def\N{\mathbb{N}}

\begin{document}

\title[Some remarks on Morse theory for posets]{Some remarks on Morse theory for posets, homological Morse theory and finite manifolds}
\author[E.G. Minian]{El\'\i as Gabriel Minian}

\address{Departamento  de Matem\'atica-IMAS\\
 FCEyN, Universidad de Buenos Aires\\ Buenos
Aires, Argentina.}

\email{gminian@dm.uba.ar}

\begin{abstract}
We introduce a version of discrete Morse theory for posets. This theory studies the topology of the order complexes $\k(X)$ of \it h-regular \rm posets $X$ from the critical points of \it admissible \rm matchings on $X$. Our approach is related to R. Forman's discrete Morse theory for CW-complexes and generalizes Forman and Chari's results on the face posets of regular CW-complexes. We also introduce a homological variant of the theory that can be used to study the topology of triangulable homology manifolds by means of their order triangulations. 
\end{abstract}




\maketitle

\section{Introduction}

Morse theory is a powerful tool to study the topology of differentiable manifolds from the critical points of smooth real-valued functions. The homotopy type of a compact manifold $N$ can be expressed in terms of the critical points of a Morse function defined on $N$. Concretely, $N$ has the homotopy type of a CW-complex with one cell of dimension $p$ for each critical point of index $p$. There are many applications of this theory in geometry, topology and physics, including its applications to surgery and cobordism, Bott periodicity theorem \cite{bott}, the h-cobordism theorem  \cite{mil2} and  J. Milnor's construction of exotic spheres \cite{mil0}. A classical reference for Morse theory is Milnor's book \cite{mil1}.

In \cite{for} R. Forman introduced a combinatorial version of Morse theory that can be applied to study the topology of polyhedra by means of discrete analogues of Morse functions and gradient vector fields defined on their triangulations. Forman's theory has, as well, a large number of applications in topology and combinatorics. M. Chari observed that Forman's discrete Morse theory for regular CW-complexes can be formulated in terms of acyclic matchings on the Hasse diagrams of their face posets \cite{cha}.  More recently  Rietsch and Williams described Forman's discrete Morse theory for general CW-complexes in terms of a certain kind of matchings on the associated posets  \cite{rw}, but the condition required is that only the regular cells are allowed to be matched. In Chari's and Rietsch-Williams' approaches, the face posets are used to restate Forman's discrete Morse theory for CW-complexes in terms of matchings and the proofs rely precisely on the fact that the posets considered are the face posets of some complexes. More concretely, the edges of the matching correspond to elementary collapses of regular faces of the CW-complex.

In this paper we study Morse theory for a  class of posets, called \it h-regular \rm posets, which includes the face posets of regular CW-complexes, and extend the main results of discrete Morse theory to this wider class. In contrast to the previous approaches, we don't require the posets to be the face posets of some complex. The Morse theory for h-regular posets study the topology of the order complexes of the posets from the critical points of the matchings. Concretely, given an h-regular poset $X$ together with an \it admissible \rm  Morse matching, we show that its order complex  $\k(X)$ is homotopy equivalent to a CW-complex with exactly one $p$-cell for each critical point of $X$ of height $p$. If $K$ is a regular CW-complex, the order complex of its face poset $\k(\x(K))$ coincides with its barycentric subdivision, therefore our approach extends Forman and Chari's result on regular CW-complexes. On the other hand, this approach provides a new insight into discrete Morse theory.  We exhibit various examples, which are not covered by the previous works on this subject, to illustrate the applications of this theory. 

In the second part of the article, we  develop a homological version of Morse theory for posets. The homological version can be applied to an even larger class, namely the class of \it cellular \rm posets. This class includes for instance the posets whose order complexes are closed homology manifolds. Given a homologically admissible Morse matching on a cellular poset $X$, one can compute the homology of $X$ with a chain complex $(\tilde C, \tilde d)$ where $\tilde C_p$ is the free abelian group generated by the critical points of $X$ of degree $p$. From this result, one can immediately deduce the well known strong and weak Morse inequalities. 

In the last section we show  examples and applications of this theory. We prove that our results can be applied to the class of h-regular CW-complexes, introduced in \cite{bm2}, without requiring regularity of the cells to be matched. Finally we show that the homological version of the theory  can be used to study homology manifolds via their order triangulations.

It is known that a finite poset is essentially the same as a finite $T_0$-topological space (see for example \cite{b1,bm1,bm2,may,mcc,sto}), and therefore this theory could also be applied to investigate the topology of finite spaces. In this article we don't adopt explicitly the finite space point of view. Moreover, the homology and the homotopy of the posets are formulated exclusively in terms of the homology and homotopy of their associated order complexes. Nevertheless  some of the constructions that appear in this article are based on results of the theory of finite spaces developed in previous papers in collaboration with J. Barmak \cite{bm1,bm2}.  

For a comprehensive exposition on discrete Morse Theory and applications,  the reader may consult Forman's articles \cite{for,for2, for3}, Chari's paper \cite{cha}, and D. Kozlov's book \cite{koz}. 
 
In this paper, all the posets, simplicial complexes and CW-complexes that we deal with, are assumed to be finite, and the homology groups $H_*(X)$ are reduced and with integer coefficients.

\section{Morse theory for posets}

We start by fixing some notation and terminology. Given a poset $X$ and an element $x\in X$,  we denote
$$\u x=\{y\in X, \  y\le x\}, \ \ut x =\set{y\in X,\  y<x}$$
and
$$\f x=\set{y\in X,\ x\leq y},\  \ft x=\set{y \in X, \ x <y}.$$

We write $x\prec y$ if $x$ is covered by $y$, i.e. if $x<y$ and there is no $z\in X$ such that $x<z<y$. The Hasse diagram $\h(X)$ of a poset $X$ is the digraph whose vertices are the points of $X$ and whose edges are the pairs $(x,y)$ such that $x\prec y$. In the graphical representation of $\h(X)$, instead of writing  the edge $(x,y)$ with an arrow, we simply put $y$ over $x$ (see Figure \ref{hreg}).

The \it order complex \rm $\k(X)$ associated to a poset $X$ is the simplicial complex of non-empty chains of $X$. This construction is functorial. It is known that  a finite poset is equivalent to a finite $T_0$-topological space, where the downsets correspond to open sets. A result of McCord \cite{mcc} asserts that a poset $X$, viewed as a finite topological space, is weak homotopy equivalent to its order complex $\k(X)$. In particular, its singular homology coincides with the homology of $\k(X)$ (cf. \cite{b1,bm1,bm2}). As we pointed out above, we will not use explicitly the finite space viewpoint and, for the purpose of this paper,  the reader may regard the homology and homotopy groups of $X$ as those of the complex $\k(X)$. 

We say that a poset $X$ is a \it model \rm of a CW-complex $Y$ if its order complex $\k(X)$ is homotopy equivalent to $Y$. The underlying polyhedron or realization of a simplicial complex $K$ will also be denoted by $K$. We say that $X$ or $\k(X)$ is an \it order triangulation \rm of its realization. Given a finite regular CW-complex $L$, the face poset $\x(L)$ is the poset of cells of $L$ ordered by inclusion. Note that $\k(\x(L))=L'$,  the barycentric subdivision of $L$.

The height $h(X)$ of a poset $X$ is the maximum of the lengths of the chains in $X$, where the length of a chain $x_0<x_1<\ldots<x_n$ is $n$. The height $h(x)$ of an element $x\in X$ is the height of $U_x$. We say that a poset $X$ is \it homogeneous \rm of dimension (or degree) $n$ if all maximal chains in $X$ have length $n$. Note that $X$ is homogeneous of dimension $n$ if and only if $\k(X)$ is a homogeneous (=pure) simplicial complex of dimension $n$. Also, a simplicial complex $L$ is homogeneous of dimension $n$ if and only if $\x(L)$ is.
A poset is \it graded \rm if $\u x$ is homogeneous for every $x\in X$. In that case, the degree of $x$, denoted by $deg(x)$,  is its height.

Note that for any regular CW-complex $L$, the poset $\x(L)$ is always graded and the degree of $x\in X$ is the dimension of $x$ (as a cell of $L$). If $X$ is homogeneous then it is graded. Moreover, a connected poset $X$ is homogeneous if and only if $X$ and its opposite $X^{op}$ are graded.

The \it join \rm $X\circledast Y$ of two posets $X$ and $Y$ is the disjoint union $X\sqcup Y$ keeping the giving ordering within $X$ and $Y$ and setting $x\leq y$ for every $x\in X$ and $y\in Y$. 
Note that $\k(X\circledast Y)=\k(X)* \k(Y)$, the simplicial join of the order complexes.

\begin{definition}\label{link1}
Let $X$ be a poset and let $x\in X$. The \it link \rm of $x$ in $X$ is the subposet $\ct x=\ut x \circledast \ft x=\set{y\in X,\ x<y \text{ or } y<x}$.
\end{definition}

\begin{remark}\label{link2}
$\k(\ct x)=lk(x,\k(X))$, the link of the vertex $x$ in $\k(X)$.
\end{remark}


\begin{definition}
A poset $X$ is called \it h-regular \rm if for every $x\in X$, $\ut X$ is a model of $S^{n-1}$ (i.e. its order complex is homotopy equivalent to $S^{n-1}$), where $n$ is the height of $x$.
\end{definition}

Note that the face poset $\x(K)$ of any regular CW-complex $K$ is h-regular. In section \ref{apps} we will deal with the face posets of h-regular CW-complexes. By definition, these posets are also h-regular.

\begin{ej} Figure \ref{hreg} shows the Hasse diagram of an h-regular poset which is not the face poset of a regular CW-complex. Note that an h-regular poset is not necessarily graded.

\begin{figure}[h]
 \begin{displaymath}
 \xymatrix@C=10pt{
& \bullet \ar@{-}[dl] \ar@{-}[dr] \ar@{-}[ddrrrr]& & & \bullet \ar@{-}[dll] \ar@{-}[dllll] \ar@{-}[ddr] & \\
\bullet \ar@{-}[d]\ar@{-}[drr]&& \bullet \ar@{-}[dll] \ar@{-}[d]
&&\\
\bullet \ar@{-}[d]\ar@{-}[drr]&& \bullet \ar@{-}[dll] \ar@{-}[d]
& & & \bullet \ar@{-}[dlll]\ar@{-}[d]\\
\bullet & & \bullet&& & \bullet 
}
\end{displaymath}
\caption{An h-regular (and non graded) poset.\label{hreg}}
\end{figure}
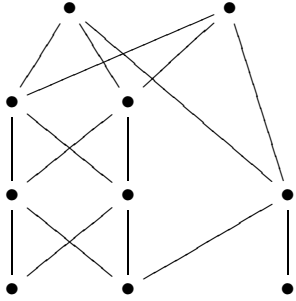
\end{ej}

\begin{definition}
We say that an edge $(x,y)\in \h(X)$ of a poset $X$ is \it admissible (in the homotopical sense) \rm  if the subposet $\ut x-\{y\}$ is homotopically trivial (i.e. its order complex is contractible). A poset is \it admissible \rm if all its edges are admissible.
\end{definition}

\begin{remark}
It is not difficult to prove that any admissible poset is automatically h-regular. Note also that the face poset $\x(K)$ of any regular CW-complex $K$ is admissible.
\end{remark}

The next example shows an admissible poset which is not the face poset of a regular CW-complex. Recall that $X=\x(K)$ for some regular CW-complex $K$ if and only if $\k(\ut x)$ is homeomorphic to a sphere for all $x\in X$ (cf. \cite[Prop 4.7.23]{blvswz}).

\begin{ej}
Let $L$ be a triangulation of the Poincar\'e homology 3-sphere and let $M=\Sigma L$ the simplicial suspension of $L$. It is known that $M$ is a 4-homology manifold which is homotopy equivalent to $S^4$ but it is not a topological manifold. Let $X=\x(M)$. Note that if $x\in X$ is a maximal point, $X-\set{x}$ is homotopically trivial since $\k(X-\set{x})$ is homotopy equivalent to $\k(X)-\set{x}$, which is the barycentric subdivision of $M$ with a point removed. This implies that the cone $Y=X\circledast *$ is admissible. Note that $Y$ is not of the form $\x(K)$ since $\k(\ut *)=\k(X)=M'$, which is not homeomorphic to a sphere.
\end{ej}

\begin{lemma}\label{altura}
Let $X$ be an h-regular poset. If an edge $(x,y)\in\h(X)$ is admissible, then $h(x)=h(y)-1$. In particular, every admissible poset is graded.
\end{lemma}

\begin{proof}
Let $h(y)=n$, then $\k(\ut y)$ is homotopy equivalent to $S^{n-1}$ and the dimension of $\k(\ut y)$ is $n-1$. Since $x\prec y$, $h(x)\leq n-1$. Suppose $h(x)< n-1$. Then any simplex in $\k(\ut y)$ which contains $x$ has dimension less than $n-1$. Therefore $\k(\ut y -x)$ has the same $(n-1)$-simplices than $\k(\ut y)$, and since $H_{n-1}(\k(\ut y))=\Z$, it follows that $H_{n-1}(\k(\ut y -x))\neq 0$, which contradicts the fact that $\k(\ut y -x)$ is contractible.
\end{proof}

Let $X$ be an h-regular poset and let $\h(X)$ be its Hasse diagram. Following Chari's construction \cite{cha}, given a matching $M$ on $\h(X)$,  let $\h_M(X)$ be the directed graph obtained from $\h(X)$ by reversing the orientations of the edges which are not in $M$. We say that a matching is a \it Morse matching \rm provided that the directed graph $\h_M(X)$ is acyclic. As usual, the points of $\h(X)$ not incident to any edge in $M$ will be called \it critical. \rm The set of critical points will be denoted by $C_M$. A Morse matching is called \it admissible \rm if every edge in $M$ is admissible. Note that, by definition, if $X$ is admissible, any Morse matching on $\h(X)$ is admissible.

\begin{lemma}\label{paths}
Let $M$ be an admissible Morse matching on $\h(X)$ and let  $x_0,\ldots,x_r$ be a directed path in $\h_M(X)$. Then $h(x_r)\leq h(x_0)+1$.
\end{lemma}
\begin{proof}
It follows from Lemma \ref{altura} and the fact that, if $x<y$, $h(x)< h(y)$.
\end{proof}

The next two lemmas follow immediately from the Gluing theorem \cite[7.5.7]{bro}. A  proof of the first lemma can be found in \cite{b2} and a proof of the second one can be found in \cite{bm2}.

\begin{lemma}\label{susp}
Let $K$ and $L$ be contractible simplicial complexes and let $T=K\cup L$. Then $T$ is homotopy equivalent to $\Sigma(K\cap L)$, the suspension of their  intersection.
\end{lemma}

\begin{lemma} \label{gamma}
Let $X$ be a poset and let $x\in X$ such that $\ct x$ is homotopically trivial, then the inclusion $\k(X - \{x\} )\hookrightarrow \k( X) $ is a homotopy equivalence.
\end{lemma}

A point $x\in X$ such that $\ct x$ is homotopically trivial, as in the preceding lemma, is called a \it $\gamma$-point \rm of $X$. This kind of points is used in \cite{bm2} to investigate reduction methods of finite topological spaces.




The following result extends the main result on discrete Morse theory to the class of h-regular posets. The techniques that we use in our proof are similar to those of Chari's \cite{cha}. However in our case the poset is not the face poset of some complex so we cannot use that the edges of the matching correspond to collapses of regular cells (cf. \cite{cha,for}). Moreover, in our context the posets are not necessarily graded, and Lemma \ref{paths} of above  is used to circumvent this problem. 

\begin{thm}\label{main}
Let $X$ be an h-regular poset and let $M$ be an admissible Morse matching on its Hasse diagram $\h(X)$. Then the order complex $\k(X)$ is homotopy equivalent to a CW-complex with exactly one $p$-cell for each critical point of $X$ of height $p$.
\end{thm}

\begin{proof}
Let $x\in X$ be a source  node of $\h_M(X)$ of maximum height. If $x$ is a maximal point of $X$, then it is critical. Let $h(x)=n$. Since $\ut x$ is a model of $S^{n-1}$ and $\k(X)=\k(X-\set{x})\cup\k(\u x)$ with $\k(X-\set{x})\cap\k(\u x)=\k(\ut x)$, by the Gluing theorem \cite[7.5.7]{bro} it turns out that $\k(X)$ is homotopy equivalent to $\k(X-\set{x})\cup e^n$, a CW-complex obtained from $\k(X-\set{x})$ by attaching an $n$-cell. On the other hand, $X-\set{x}$ is still h-regular and the admissible Morse matching $M$ restricts to an admissible Morse matching on $\h(X-\set{x})$. 

If $x$ is not a maximal point, there exists $y\in X$ such that $x\prec y$. Since $x$ is a source node, it follows that $(x,y)\in M$ and that there is no $z\neq y$ such that $x\prec z$. Suppose there is an element $w\in X$ such that $x\prec y\prec w$. Since $(x,y)\in M$, then $(y,w)\notin M$. Let $s$ be a source node such that $s=x_0,\ldots,x_r=w$ is a directed path in $\h_M(X)$. Then, by Lemma \ref{paths} $h(s)\geq h(w)-1\geq h(x)+1$, which contradicts the maximality of the height of $x$. It follows that $\ft x=\set{y}$ with $y$ a maximal element of $X$. Note that the link $lk(x,\k(X))$ is a cone. In particular, $\ct x$ is homotopically trivial and, by Lemma \ref{gamma},  $\k(X-\set{x})\subset \k(X)$ is a strong deformation retract.

Note also that $\ct y^{X-\set{x}}=\ut y^{X-\set{x}}=\ut y^X-\set{x}$, which is homotopically trivial. Here we write $\ct y^{X-\set{x}}, \ \ut y^{X-\set{x}}$ and $\ut y^X$ to distinguish whether the subposets $\ct y, \ut y$ are considered in $X-\set{x}$ or in $X$.

 Again  by Lemma \ref{gamma} we have
$$\k(X-\set{x,y})\simeq \k(X-\set{x})\simeq \k(X).$$
Therefore we can remove $x$ and $y$ from $X$ without affecting the homotopy type of the order complex. On the other hand, since $y$ is a maximal point of $X$ and $\ft x=\set{y}$, the subposet $X-\set{x,y}$ is still h-regular and the admissible Morse matching $M$ on $\h(X)$ restricts to an admissible Morse matching on $\h(X-\set{x,y})$.

The result now follows inductively.
\end{proof}

\begin{remark}
Note that Theorem \ref{main} extends the theory of Forman and Chari for regular CW-complexes since the face posets $\x(K)$ of regular CW-complexes are admissible and the order complex $\k(\x(K))$ of a regular CW-complex $K$ is the barycentric subdivision of $K$.
\end{remark}

A point $x\in X$ such that $\ft x=\set{y}$ (or more generally, if $\ft x$ has a minimum), as in the proof of Theorem \ref{main}, is called an \it up beat point \rm or a \it linear point \rm  of $X$. This kind of points was studied by R. Stong \cite{sto} to classify the homotopy types of finite spaces (see also \cite{b1,bm1,bm2,may}).

We finish this section with two examples of h-regular posets with admissible Morse matchings, to illustrate the situation. Note that they are not the face posets of regular complexes.

\begin{ej} Figure \ref{admatch} shows an admissible Morse matching on the h-regular poset $X$ of figure \ref{hreg}. The edges of the matching are represented with dashed arrows. Its order complex is homotopy equivalent to a CW-complex with one $0$-cell and one $3$-cell, thus $X$ is a model of $S^3$.

\begin{figure}[h]
 \begin{displaymath}
 \xymatrix@C=10pt{
& \bullet  \ar@{-}[dr] \ar@{-}[ddrrrr]& & & \bullet \ar@{-}[dll] \ar@{-}[dllll] \ar@{-}[ddr] & \\
\bullet \ar@{-->}[ur] \ar@{-}[d]\ar@{-}[drr]&& \bullet  \ar@{-}[d]
&&\\
\bullet \ar@{-->}[urr] \ar@{-}[d]\ar@{-}[drr]&& \bullet  \ar@{-}[d]
& & & \bullet \ar@{-}[d]\\
\bullet\ar@{-->}[urr] & & \bullet \ar@{-->}[urrr] && & \bullet 
}
\end{displaymath}
\caption{An admissible Morse matching on an h-regular model of $S^3$.\label{admatch}}
\end{figure}
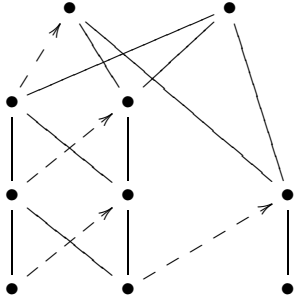

\end{ej}

\begin{ej} The admissible Morse matching on the poset of Figure \ref{contrac} shows that its order complex is contractible.

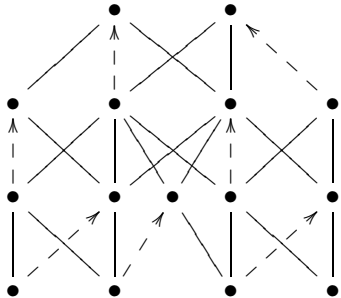
\begin{figure}[h]
 \begin{displaymath}
 \xymatrix@C=10pt{
& & \bullet  \ar@{-}[drr] \ar@{-}[dll]& &  \bullet \ar@{-}[dll] \ar@{-}[d]  & & \\
\bullet  \ar@{-}[drr]&& \bullet  \ar@{-}[d]\ar@{-->}[u]\ar@{-}[dll]\ar@{-}[dr]\ar@{-}[drr]
&& \bullet \ar@{-}[dll]\ar@{-}[dl]\ar@{-}[drr]&& \bullet\ar@{-}[dll]\ar@{-}[d]\ar@{-->}[ull]\\
\bullet \ar@{-->}[u] \ar@{-}[d]\ar@{-}[drr]&& \bullet  \ar@{-}[d]
& \bullet \ar@{-}[dr]& \bullet\ar@{-->}[u] \ar@{-}[d]\ar@{-}[drr]& & \bullet \ar@{-}[d]\\
\bullet\ar@{-->}[urr] & & \bullet \ar@{-->}[ur] &&  \bullet \ar@{-->}[urr]&& \bullet
}
\end{displaymath}
\caption{An admissible Morse matching on a homotopically trivial h-regular poset.\label{contrac}}
\end{figure}
\end{ej}

\section{Homological Morse theory and cellular posets}

In this section we will investigate a homological variant of our theory. First we introduce the notion of \it cellular \rm poset and show that the homology of these posets can be computed using the \it cellular chain complex \rm $(C_*,d)$ which, in degree $p$, consists of the free abelian group generated by the elements of the poset of degree $p$. In the second part of this section we construct the Morse complex of a cellular poset with a homologically admissible matching, adapting Forman's techniques to this context \cite[Prop 6.3 to Thm 8.2]{for}.

\begin{definition}
A  \it cellular \rm poset is a graded poset $X$ such that for every $x\in X$, $\ut x$ has the homology of a $(p-1)$-sphere, where $p=deg(x)$.
\end{definition}

\begin{remark}
Note that, by definition, any graded h-regular poset is cellular.
\end{remark}

The \it $p$-skeleton \rm of a graded poset $X$ is the subposet $X^p=\set{ x\in X,\ deg(x)\leq p}$. There is a filtration
$$X^0\subseteq X^1\subseteq\ldots\subseteq X^n=X.$$

Note that if $X=\x(L)$ is the poset of cells of a regular CW-complex $L$, then $X^p=\x(L^p)$, the poset associated to the $p$-skeleton of $L$.

The proof of the following result is routine and therefore omitted.

\begin{proposition}
Let $X$ be a graded poset and let $X^p$ be its $p$-skeleton. Then
$$H_r(X^p,X^{p-1})=\bigoplus_{deg(x)=p}H_r(\u x,\ut x)=\bigoplus_{deg(x)=p}H_{r-1}(\ut x).$$
\end{proposition}

Applying the long exact sequence of the pair $(X^{p+1},X^p)$, we deduce the following
 
\begin{corollary}\label{homologiabaja}
Let $X$ be a cellular poset. Then, the inclusion $X^p\subseteq X$ induces an isomorphism $H_r(X^p)=H_r(X)$ for any $r<p$.
\end{corollary}






\begin{remark}\label{homologiaalta}
Note that, by a dimension argument, if $X$ is a graded poset and $r>p$, then $H_r(X^p)=0$.
\end{remark}

\begin{definition}
Given a cellular poset $X$, its \it cellular chain complex \rm  $(C_*,d)$ is defined as follows
$$C_p:=H_p(X^p,X^{p-1})=\bigoplus_{deg(x)=p}H_{p-1}(\ut x)$$
which is a free abelian group with one generator for each element of $X$ of degree $p$. The differential $d:C_p\to C_{p-1}$ is defined as the composition
$$\xymatrix{H_{p}(X^{p},X^{p-1})\ar[r]^-{\partial} & H_{p-1}(X^{p-1}) \ar[r]^-j & H_{p-1}(X^{p-1},X^{p-2})}.$$
Here $j$ is the canonical map induced by the inclusion.
\end{definition}

\begin{thm}
Let $X$ be a cellular poset and let $(C_*,d)$ be its cellular chain complex. Then $H_*(C_*)=H_*(X)$.
\end{thm}
\begin{proof}
Follows  from \ref{homologiabaja} and \ref{homologiaalta}.
\end{proof}

We choose a generator of $H_{p-1}(\ut x)=\Z$ for every $x\in X$ of degree $p$ and then we identify $C_p$ with the free abelian group with basis $\set{x\in X, deg(x)=p}$. It is not hard to prove that the differential $d:C_p\to C_{p-1}$ has the form $$d(x)=\displaystyle\sum_{w\prec x}\epsilon(x,w)w$$ where the incidence number $\epsilon(x,w)\in \Z$ is the degree of the map
$$\tilde\partial:\Z =H_{p-1}(\ut x)\to H_{p-2}(\ut w)=\Z,$$
which coincides with the connecting morphism of the Mayer-Vietoris sequence associated to the covering $\ut x=(\ut x-\set{w})\cup \u w$.
This means that $\tilde\partial(x)=\epsilon(x,w)w$, where $x$ and $w$ represent the chosen generators of $H_{p-1}(\ut x)$ and $H_{p-2}(\ut w)$ respectively.


\begin{definition}
Let $X$ be a poset. An edge $(w,x)\in \h(X)$ is  \it homologically admissible \rm  if $\ut x-\set{w}$ is acyclic. 
\end{definition}

A poset is homologically admissible if all its edges are homologically admissible. It is not difficult to see that such posets are cellular.

\begin{remark}\label{iscellular}
If $(w,x)$ is a homologically admissible edge of a cellular poset $X$, the incidence number $\epsilon(x,w)$ is $1$ or $-1$. This follows from the fact that the map $\tilde\partial$ of above is an isomorphism, since $H_*(\ut x-\set{w})=H_*(\u w)=0$. 
\end{remark}

\begin{remark}
Note that an admissible poset in the homotopical sense is homologically admissible. An example of a homologically admissible poset which is not admissible in the homotopical sense is the cone $\x(L)\circledast *$ of the poset $\x(L)$ associated to a triangulation of the Poincar\'e homology 3-sphere $L$ (see subsection \ref{manifolds}).
\end{remark}

We will show that the homology of a cellular poset $X$ can be computed, from a homologically admissible Morse matching $M$,  with a chain complex $(\tilde C_*,\tilde d)$ which in degree $p$ consists of the abelian group generated by the critical points of $X$ of degree $p$. This extends the classical and Forman's theories to the class of cellular posets. In particular, we will deduce the strong and weak Morse inequalities similarly as in the classical theory.

Our construction of the Morse chain complex will follow Forman's \cite{for}. We need first to define the \it Morse function \rm $f:X\to\R$ associated to a Morse matching $M$. This function will be used as an auxiliary tool in the proofs. 

\begin{definition} A Morse function $f:X\to\R$ is a set theoretic function such that for every $x\in X$
\begin{itemize}
\item $\#\set{y\in X,\ x\prec y \text{ and } f(x)\geq f(y)}\leq 1$ and
\item $\#\set{z\in X,\ z\prec x \text{ and } f(z)\geq f(x)}\leq 1$.
\end{itemize}
A point $x\in X$ is critical if 
\begin{itemize}
\item $\#\set{y\in X,\ x\prec y \text{ and } f(x)\geq f(y)}=0$ and
\item $\#\set{z\in X,\ z\prec x \text{ and } f(z)\geq f(x)}=0$.
\end{itemize}
The set of critical points of $f$ is denoted by $C_f$.
\end{definition}

Given a Morse matching $M$ on $\h(X)$, a \it Morse path \rm  is a path in $\h_M(X)$ of the form $x_0,y_0,x_1,y_1,\ldots,x_r$ such that $(x_i,y_i)\in M$. Note that $x_{i+1}\prec y_i$ for every $i$. The length of the path $x_0,y_0,x_1,y_1,\ldots,x_r$ is $r$. We define $l_M(x)\in \N_0$ to be the maximum of the lengths of the Morse paths which start in $x$. Note that, if $x$ is a critical point,  $l_M(x)=0$.

\begin{lemma}\label{funct}
Given a Morse matching on a cellular poset $X$, there is a Morse function $f$ on $X$ such that $C_f=C_M$.
\end{lemma}
\begin{proof}
The proof is similar to Forman's \cite{for} (see also \cite{cap} for an alternative proof). We construct by induction a function $f_r:X^r\to \R$, where $X^r$ is the $r$-skeleton of $X$. Let $f_0:X^0\to \R$ be the constant map $0$. Suppose we have already defined $f_{r-1}:X^{r-1}\to\R$. Let $L=max \set{l_M(x), deg(x)=r-1}$. Define $f_r$ as follows.

$$f_r(x)=\begin{cases}
f_{r-1}(x) & deg(x)\leq r-2\\
f_{r-1}(x)+\frac{l_M(x)}{L+1} & deg(x)=r-1\\
r & deg(x)=r \text{ and } \nexists w\prec x \text{ such that } (w,x)\in M\\
f_{r-1}(w)+\frac{l_M(w)}{L+1} & deg(x)=r \text{ and } (w,x)\in M
\end{cases}$$
\end{proof}

\begin{ej} Figure \ref{function} shows the Morse function associated to a Morse matching.

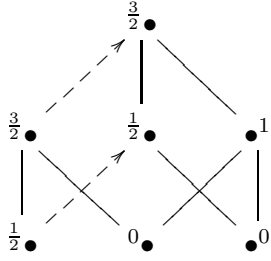
\begin{figure}[h]
 \begin{displaymath}
 \xymatrix@C=10pt{
& & ^{\frac{3}{2}} \bullet  \ar@{-}[d] \ar@{-}[drr]\\
 ^{\frac{3}{2}} \bullet  \ar@{-}[d] \ar@{-}[drr] \ar@{-->}[urr] & &   ^{\frac{1}{2}} \bullet  \ar@{-}[drr] & & \bullet ^1 \ar@{-}[d] \ar@{-}[dll] \\
  ^{\frac{1}{2}} \bullet  \ar@{-->}[urr] && ^0 \bullet && \bullet ^0
}
\end{displaymath}
\caption{The Morse function associated to a Morse matching.\label{function}}
\end{figure}

\end{ej}

We define now the \it Morse complex \rm associated to a cellular poset $X$ with a homologically admissible Morse matching $M$ following Forman's construction \cite[Prop 6.3 to Thm 8.2]{for}. In order to adapt Forman's construction and proofs to our setting, we use Remark \ref{iscellular} and the Morse function $f$ associated to $M$, constructed in Lemma \ref{funct}.

Following Forman's notation, we define a map $V:C_p\to C_{p+1}$ as follows.
$$V(x)=\begin{cases}
-\epsilon(y,x) y & \text{ if there is  } y\in X \text{ with } (x,y)\in M\\
0 &  \text{ otherwise.} \end{cases}$$

The \it gradient flow \rm $\phi:C_*\to C_*$ is defined as $\phi=1+dV+Vd$.
Since $d\phi=\phi d$, the $\phi$-invariant chains $C_*^{\phi}$ form a well defined subcomplex. Concretely, 
$$C_p^{\phi}(X)=\set{c\in C_p(X),\ \phi(c)=c}.$$

The complex $(C_*^{\phi},d)$ is the \it Morse complex \rm  of the cellular poset $X$.

Note that Theorems 6.4, 7.1, 7.2 and 7.3 of \cite{for} remain valid in our context. The proofs are exactly the same as in the case of regular CW-complexes. In our setting we use the cellular chain complex of the poset and the associated Morse function $f$ of Lemma \ref{funct}. Therefore we have proved

\begin{thm}
If $X$ is a cellular poset and $M$ is a homologically admissible Morse matching on $\h(X)$, the homology of the associated Morse complex $(C_*^{\phi},d)$ coincides with the homology of $X$.
\end{thm}

On the other hand, similarly as in Forman's setting, for each $p$, the abelian group 
$C_p^{\phi}$ is isomorphic to the free abelian group $\tilde C_p$ spanned by the critical $p$-cells. This follows from the results 7.1, 8.1 and 8.2 of \cite{for}. The extension of Forman's theory to cellular posets follows immediately.

\begin{corollary}
If $X$ is a cellular poset with a homologically admissible Morse matching $M$, the homology of $X$ can be computed with a chain complex $(\tilde C, \tilde d)$, where $\tilde C_p$ is the free abelian group generated by the critical points of $X$ of degree $p$. In particular, the strong Morse inequalities (and hence, the weak Morse inequalities) hold.
\end{corollary}

\section{Examples and applications}\label{apps}

\subsection{h-regular CW-complexes}

The concept of \it h-regular \rm CW-complex was introduced in \cite{bm2} (see also \cite{b1}). It generalizes the notion of a regular complex and allows one to manipulate less rigid CW-structures than the regular ones (and therefore with fewer cells), but with similar nice properties as the properties of the regular complexes.

A CW-complex $K$ is \it h-regular \rm if the attaching map of each cell $e^n$ is a homotopy equivalence onto its image $\dot{e}^n$ and the closed cells $\overline{e^n}$ are subcomplexes of $K$. Equivalently, $K$ is h-regular if the closed cells are contractible subcomplexes. The reader can find various examples of h-regular structures on CW-complexes in \cite{b1,bm2}.

If $K$ is a regular complex, $\k(\x(K))=K'$, the barycentric subdivision of $K$. If $K$ is any CW-complex, $\k(\x(K))$ is not in general homeomorphic to $K$, they are not even homotopy equivalent. The following result, which was proved in \cite{bm2}, allows us to apply the results of the previous sections to the h-regular CW-complexes.
\begin{proposition}
If $K$ is an h-regular CW-complex, $\k(\x(K))$ and $K$ are homotopy equivalent. 
\end{proposition}

The last proposition asserts that the face poset of an h-regular complex is a model of the complex. Since, by definition, the face poset of an h-regular CW-complex is h-regular, then by 
Theorem \ref{main}, we obtain the following result.

\begin{corollary}\label{main2}
Let $K$ be an h-regular CW-complex and let $M$ be an admissible Morse matching on its face poset $\x(K)$. Then $K$ is homotopy equivalent to a CW-complex with exactly one $p$-cell for each critical $p$-cell of $K$.
\end{corollary}

The next example illustrates this situation.

\begin{ej}
Consider the following h-regular structure on the space $K$ which is obtained from the square by identifying all its boundary edges as indicated in the picture.

\begin{center}
\includegraphics[scale=0.5]{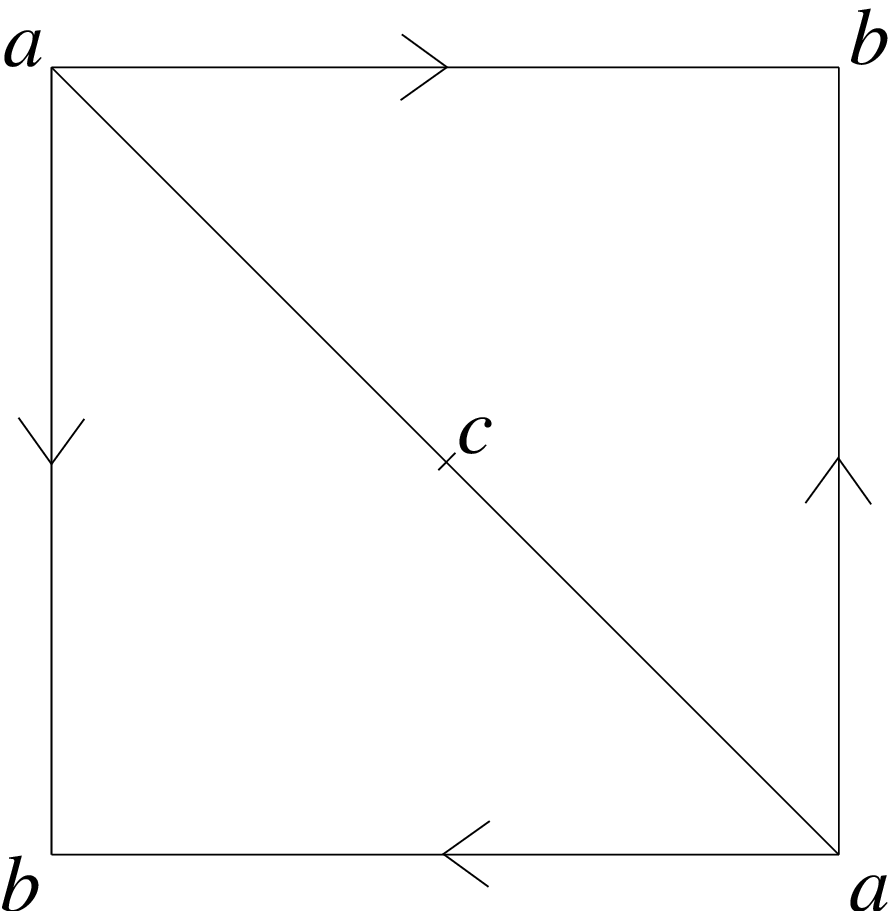}
\end{center}

The following admissible Morse matching on $\x(K)$ shows that $K$ is homotopy equivalent to $S^2$ (cf. \cite[Example 4.8]{bm2}).

\begin{displaymath}
 \xymatrix@C=10pt{
\bullet \ar@{-}[drrrr] \ar@{-}[d]\ar@{-}[drr]&& \bullet  \ar@{-}[d]  \ar@{-}[drr]
& &\\
\bullet \ar@{-->}[urr] \ar@{-}[d]\ar@{-}[drr]&& \bullet  \ar@{-}[d]
& &  \bullet ^{ab}\ar@{-}[d]\\
_c \bullet\ar@{-->}[urr] & & \bullet _a \ar@{-->}[urr] & & \bullet _b
}
\end{displaymath}

\end{ej}
 
 Note that Corollary \ref{main2} cannot be deduced from  Forman  and Chari's results, nor from Rietsch-Williams's approach. In our context we allow non-regular cells to be matched.
 
\subsection{Finite manifolds and finite homology manifolds }\label{manifolds}

\bigskip

An interesting class of admissible posets is constituted by the \it closed finite manifolds. \rm A finite manifold is a poset $X$ such that its order complex $\k(X)$ is a combinatorial manifold. An alternative definition is the following.

\begin{definition}
A poset $X$ is called an \it $n$-sphere \rm if its order complex $\k(X)$ is a combinatorial $n$-sphere, i.e. if $\k(X)$ is PL-homeomorphic to the boundary of an $(n+1)$-simplex. A poset $X$ is  an \it $n$-ball \rm  if $\k(X)$ is a combinatorial $n$-ball, i.e. if it is PL-homeomorphic to an $n$-simplex. A poset $X$ is a \it finite manifold \rm  of dimension $n$, or an \it $n$-manifold, \rm   if it is homogeneous of dimension $n$ and for each $x\in X$ of degree $0$, $\ft x$ is a ball or a sphere of dimension $n-1$ and for each $x\in X$ of degree $n$, $\ut x$ is a ball or a sphere of dimension $n-1$.
\end{definition}

Recall that a combinatorial $n$-manifold is a simplicial complex $M$ such that for every vertex $x\in M$, $lk(x,M)$ is a sphere or a ball of dimension $n-1$. For a comprehensive exposition of the theory of combinatorial manifolds, we refer the reader to \cite{lic}.

 The following proposition relates finite manifolds with combinatorial manifolds (cf. \cite{nt}).

\begin{proposition}
A poset $X$ is an $n$-manifold if and only if $\k(X)$ is a combinatorial $n$-manifold.
\end{proposition}
\begin{proof}
Since the join of a ball or a sphere with a ball or a sphere is again a ball or sphere, by Remark \ref{link2} it suffices to prove that if $X$ is an $n$-manifold and $x\in X$, the subposets $\ut x$ and $\ft x$ are balls or spheres. We prove this by induction on $n$.
Suppose the assertion is true for manifolds of dimension $n-1$ and let $x\in X$, with $X$ a manifold of dimension $n$. If $x$ is a maximal or minimal point of $X$, then there is nothing to prove.
In other case, let $y\in X$ be a minimal element of $X$ which is comparable with $x$. Since $\ft y$ is ball or sphere of dimension $n-1$, in particular it is an $(n-1)$-manifold and, by induction $\ft x$ is a ball or sphere. Similarly, consider a maximal element $z\in X$ comparable with $x$.   Since $\ut z$ is ball or sphere of dimension $n-1$, in particular it is an $(n-1)$-manifold and, by induction $\ut x$ is a ball or sphere.
\end{proof}

A finite manifold $X$ is \it closed \rm if $\k(X)$ is a closed manifold, i.e. if the links of the points are spheres. If $X$ is an $n$-sphere and $x\in X$, then by Newman's Theorem, $X-\set{x}$ is an $n$-ball (cf. \cite{lic}). It follows that closed finite manifolds are admissible. In fact, by \cite[Prop 4.7.23]{blvswz}, it is easy to see that a closed manifold has the form $\x(K)$ for some regular CW-complex $K$.

We focus our attention now on the class of finite closed homology manifolds. A poset $X$ is a \it finite closed homology manifold \rm if its order complex $\k(X)$ is a closed homology manifold, i.e. if the link of every simplex $\sigma\in\k(X)$ has the homology of $S^{n-k-1}$, where $n$ is the dimension of the manifold and $k$ the dimension of the simplex. It is clear that the links of homology manifolds are also homology manifolds. 

We will show that closed homology manifolds are homologically admissible. This allows one to use the homological version of Morse theory developed in the previous section to investigate the topology of triangulable homology manifolds by means of their order triangulations. In general a closed homology manifold is not necessarily  admissible in the homotopical sense. We show below an example of an order triangulation of a homology sphere which is not admissible in the homotopical sense.

\begin{thm}
If $X$ is a closed homology manifold, it is homologically admissible. 
\end{thm}

\begin{proof}
Let $X$ be a closed homology $n$-manifold and let $x\in X$ of degree $r$. Note that $\ut x$ is a homology sphere, i.e. a homology manifold with the homology of a sphere, since its order complex is the link in $\k(X)$ of the simplex associated to any maximal chain $x=x_r<x_{r+1}<\ldots<x_n$ with minimum $x$. If $y\in\ut x$ is a maximal element, $\k(\ut x-\set{y})\subset \k(\ut x)-\set{y}$ is a strong deformation retract and, since a homology sphere with a point removed is acyclic, it follows that $X$ is homologically admissible.
\end{proof}

\begin{ej}
Let $L$ be a triangulation of the Poincar\'e homology 3-sphere. The poset $\x(L)\circledast \x(L)$ is a homology sphere, since $\k(\x(L)\circledast \x(L))=L'*L'$. On the other hand, $\x(L)\circledast \x(L)$ is not admissible in the homotopical sense, since for every minimal point $x$ of the second copy of $\x(L)$, $\ut x^{\x(L)\circledast \x(L)}=\x(L)$ which is not simply connected.

\end{ej}

{\bf Acknowledgement.}
I would like to thank Jonathan Barmak for very useful comments during the preparation and correction of this paper.


\begin{thebibliography}{xxx}

\bibitem{b1} J.A. Barmak. \textit{Algebraic topology of finite topological spaces and applications}. Lecture Notes in Mathematics Vol. 2032. Springer (2011).

\bibitem{b2} J.A. Barmak. \textit{Star clusters in independence complexes of graphs}. Preprint (2010) arXiv:1007.0418v1.

\bibitem{bm1} J.A. Barmak and E.G. Minian. \textit{Simple homotopy types and finite spaces}. Adv. Math. 218 (2008), no. 1, 87-104.


\bibitem{bm2} J.A. Barmak and E.G. Minian. \textit{One-point reductions of finite spaces, h-regular CW-complexes and collapsibility}.
		Algebr. Geom. Topol. 8 (2008), no.3, 1763-1780.


\bibitem{blvswz} A. Bj\"orner, M. Las Vergnas, B. Sturmfelds, N.White and G. Ziegler. \textit{Oriented matroids}. Enclyclopedia of Mathematics and its Applications. Cambridge University Press (1999).

\bibitem{bott}  R. Bott. \textit{The stable homotopy of the classical groups}. Ann. of Math.(2) 70 (1959), 313-337.

\bibitem{bro} R. Brown. \textit{Topology and groupoids}. BookSurge LLC (2006).
 
\bibitem{cap} N. Capitelli. \textit{Colapsabilidad en variedades combinatorias y espacios de deformaciones}. Diploma Thesis, Universidad de Buenos Aires (2009). Available at \url{http://cms.dm.uba.ar/academico/carreras/licenciatura/tesis/capitelli.pdf}.

\bibitem{cha} M. Chari. \textit{On discrete Morse functions and combinatorial decompositions}, Formal power series and algebraic combinatorics (Vienna, 1997). Discrete Math. 217 (2000), no. 1-3, 101-113.

\bibitem{for} R. Forman. \textit{Morse theory for cell complexes}. Adv. Math. 134 (1998),no.1, 90-145.

\bibitem{for2} R. Forman. \textit{Witten-Morse theory for cell complexes}. Topology 37 (1998), no. 5, 945-979.

\bibitem{for3} R. Forman. \textit{A user's guide to discrete Morse theory}. Sem. Lothar. Combin. 48 (2002), Art. B48c, 35 pp.

\bibitem{gs} D. Galewski and R. Stern. \textit{Classification of simplicial triangulations of topological manifolds}. Ann. of Math. 111 (1980) 1-34.

\bibitem{koz} D. Kozlov. {\it Combinatorial algebraic topology}. Algorithms and Computation in Mathematics, Vol. 21. Springer, Berlin (2008).

\bibitem{lic} W. Lickorish. {\it Simplicial moves on complexes and manifolds}. Geom.  Topol. Monographs. Vol 2. (1999) 299-320.

\bibitem{may} J.P. May. \textit{Finite topological spaces}.
    Notes for REU (2003).  Available at \url{http://www.math.uchicago.edu/~may/MISCMaster.html}.

\bibitem{mcc} M.C. McCord. \textit{Singular homology groups and homotopy groups of finite topological spaces}.
    Duke Math. J. 33 (1966), 465-474.
    
\bibitem{mil0}  J. Milnor. \textit{On manifolds homeomorphic to the $7$-sphere}. Ann. of Math. 64 (1956), 399-405.

\bibitem{mil1} J. Milnor. \textit{Morse theory. Based on lecture notes by M. Spivak and R. Wells}. Annals of Mathematics Studies, No. 51 Princeton University Press, Princeton, N.J. (1963) vi+153 pp.    

\bibitem{mil2} J. Milnor. \textit{Lectures on the $h$-cobordism theorem}.Princeton University Press, Princeton, N.J. (1965) v+116 pp. 
   



\bibitem{nt} S. Negami and M. Tsuchiya. {\it Manifold posets}. Sci. Rep. Yokohama Nat. Univ. Sect. I Math. Phys. Chem. No. 41 (1994), 23-32. 

\bibitem{rw} K. Rietsch and L. Williams. {\it Discrete Morse theory for totally non-negative flag varieties}.  Adv. Math. 223 (2010), no. 6, 1855-1884. 

\bibitem{sto} R.E. Stong. \textit{Finite topological spaces}.
    Trans. Amer. Math. Soc. 123 (1966), 325-340.


\end{thebibliography}
\end{document}